\documentclass[draft]{amsart}
\usepackage{amssymb}
\usepackage{braket}
\usepackage{mathrsfs}
\usepackage[all]{xy}

\newtheorem{thm}{Theorem}[section]
\newtheorem{cor}[thm]{Corollary}
\newtheorem{prop}[thm]{Proposition}
\theoremstyle{definition}
\newtheorem{dfn}[thm]{Definition}
\newtheorem{ex}[thm]{Example}
\newtheorem{claim}[thm]{Claim}
\newtheorem{fact}[thm]{Fact}
\theoremstyle{remark}
\newtheorem{rem}[thm]{Remark}

\newcommand{\Ob}{\mathrm{Ob}}
\newcommand{\id}{\mathrm{id}}
\newcommand{\ppr}{^{\prime}}

\newcommand{\Mon}{\mathit{Mon}}
\newcommand{\Sett}{\mathit{Set}}
\newcommand{\Ab}{\mathit{Ab}}
\newcommand{\Add}{\mathit{Add}}
\newcommand{\SRing}{\mathit{SRing}}
\newcommand{\Ring}{\mathit{Ring}}

\newcommand{\TamG}{\mathit{Tam}(G)}
\newcommand{\HopfG}{\mathit{Hopf}(G)}
\newcommand{\STamG}{\mathit{STam}(G)}
\newcommand{\MackG}{\mathit{Mack}(G)}
\newcommand{\SMackG}{\mathit{SMack}(G)}
\newcommand{\SMackK}{\mathit{SMack}(K)}
\newcommand{\GreenG}{\mathit{Green}(G)}
\newcommand{\Gs}{{}_G\mathit{set}}
\newcommand{\ev}{\mathit{ev}^G}

\newcommand{\GMon}{G\text{-}\mathit{Mon}}
\newcommand{\Gmon}{G\text{-}\mathit{mon}}
\newcommand{\RAlg}{R\text{-}\mathit{Alg}}

\newtheorem{thmtw}{{\rm \textbf{Theorem \ref{MainThm}}}}

\newtheorem{thmto}{{\rm \textbf{Theorem \ref{MainThm} {\rm (iv)}}}}

\newtheorem{thmlast}{{\rm \textbf{Theorem \ref{ThmSimulPoly}}}}

\newcommand{\TTamG}{T\text{-}\mathit{Tam}(G)}

\newcommand{\HH}{\mathcal{H}}
\newcommand{\M}{\mathcal{M}}
\newcommand{\U}{\mathscr{U}}
\newcommand{\V}{\mathscr{V}}

\numberwithin{equation}{section}

\begin{document}

\title[A generalization of the Dress construction]{A generalization of The Dress construction for a Tambara functor, and polynomial Tambara functors}

\author{Hiroyuki NAKAOKA}
\address{Department of Mathematics and Computer Science, Kagoshima University, 1-21-35 Korimoto, Kagoshima, 890-0065 Japan}

\email{nakaoka@sci.kagoshima-u.ac.jp}

\thanks{The author wishes to thank Professor Fumihito Oda for stimulating arguments and useful comments}
\thanks{The author wishes to thank Professor Serge Bouc and Professor Radu Stancu for stimulating arguments and their useful comments and advices}
\thanks{Supported by JSPS Grant-in-Aid for Young Scientists (B) 22740005}

\begin{abstract}
For a finite group $G$, 
(semi-)Mackey functors and (semi-)Tambara functors are regarded as $G$-bivariant analogs of (semi-)groups and (semi-)rings respectively. In fact if $G$ is trivial, they agree with the ordinary (semi-)groups and (semi-)rings, and many naive algebraic properties concerning rings and groups have been extended to these $G$-bivariant analogous notions.

In this article, we investigate a $G$-bivariant analog of the semi-group rings with coefficients. Just as a coefficient ring $R$ and a monoid $Q$ yield the semi-group ring $R[Q]$, our constrcution enables us to make a Tambara functor $T[M]$ out of a semi-Mackey functor $M$, and a coefficient Tambara functor $T$. 
This construction is a composant of the Tambarization and the Dress construction.

As expected, this construction is the one uniquely determined by the righteous adjoint property. Besides in analogy with the trivial group case, if $M$ is a Mackey functor, then $T[M]$ is equipped with a natural Hopf structure.

Moreover, as an application of the above construction, we also obtain some $G$-bivariant analogs of the polynomial rings. 
\end{abstract}

\maketitle


\section{Introduction and preliminaries}


For a finite group $G$, a (resp. semi-)Mackey functor is a pair of a contravariant functor and a covariant functor to the category of abelian groups $\Ab$ (resp. of commutative monoids $\Mon$), satisfying some conditions (Definition \ref{DefMackFtr}). Since the category of Mackey functors is a symmetric monoidal abelian category which agrees with $\Ab$ when $G$ is trivial, it is regarded as a $G$-bivariant analog of $\Ab$. Similarly a semi-Mackey functor is regarded as a $G$-bivariant analog of a commutative monoid.

In this view, a Tambara functor is regarded as a $G$-bivariant analog of a commutative ring. It consists of an additive Mackey functor structure and a multiplicative semi-Mackey functor structure, satisfying the \lq distributive law' (Definition \ref{DefTamFtr}).

Some naive algebraic properties concerning rings and groups have been extended to these $G$-bivariant analogous notions. For example, in our previous article (\cite{N_TamMack}), as a $G$-bivariant analog of the functor taking semi-group rings
\[ \mathbb{Z}[-]\colon\Mon\rightarrow\Ring, \]
we constructed a functor called the {\it Tambarization functor}
\[ \mathcal{T}\colon\SMackG\rightarrow\TamG, \]
which is characterized by some natural adjoint property (Fact \ref{FactNTM1}). Here $\SMackG$ denotes the category of semi-Mackey functors, and $\TamG$ denotes the category of Tambara functors.

In this article, more generally, we investigate a $G$-bivariant analog of the semi-group ring {\it with a coefficient ring}. In the trivial group case, from any commutative ring $R$ and any commutative monoid $Q$ we can make the semi-group ring $R[Q]$,
and this gives a functor
\[ \Ring\times\Mon\rightarrow\Ring\ ;\ (R,Q)\mapsto R[Q]. \]
In section \ref{Section2}, analogously we construct a functor
\[ \TamG\times\SMackG\rightarrow\TamG, \]
which unifies the Tambarization (\cite{N_TamMack}) and the Dress construction (\cite{O-Y3}) as follows: 
\begin{thmtw}
For any finite group $G$, there is a functor
\[ \mathcal{F}\colon\TamG\times\SMackG\rightarrow\TamG \]
which satisfies the following.
\begin{enumerate}
\item[{\rm (i)}] If $G$ is trivial, then $\mathcal{F}$ agrees with the functor taking semi-group rings with coefficients
\[ \Ring\times\Mon\rightarrow\Ring\ ;\ (R,Q)\mapsto R[Q]. \]
\item[{\rm (ii)}] If $T=\Omega$, we have a natural isomorphism $\mathcal{F}(\Omega,M)\cong\mathcal{T}(M)$ for each semi-Mackey functor $M$. 
\item[{\rm (iii)}] If $Q$ is a finite $G$-monoid, then we have a natural isomorphism $\mathcal{F}(T,\mathcal{P}_Q)\cong T_Q$ for each Tambara functor $T$. Here, $T_Q$ is the Tambara functor obtained through the Dress construction (\cite{O-Y3}).
\end{enumerate}
\end{thmtw}

As expected from the trivial group case, $\mathcal{F}$ is the unique functor characterized by the following adjoint property:
\begin{thmto}
For each Tambara functor $T$ and semi-Mackey functor $M$, naturally $\mathcal{F}(T,M)$ becomes a $T$-Tambara functor. Moreover if we fix $T$, then the induced functor
\[ \mathcal{F}(T,-)\colon\SMackG\rightarrow\TTamG \]
is left adjoint to the composition of forgetful functors
\[ \TTamG\rightarrow\TamG\overset{(-)^{\mu}}{\longrightarrow}\SMackG\ ;\ (S,\sigma)\mapsto S^{\mu}. \]
\end{thmto}

\smallskip

\noindent Here $\TTamG$ denotes the category of $T$-Tambara functors (Definition \ref{DefT-Tam}). Besides in analogy with the trivial group case, if $M$ is a Mackey functor, then $T[M]=\mathcal{F}(T,M)$ is equipped with a natural Hopf structure $(T[M],\Delta_M,\varepsilon_M,\eta_M)$ $($Corollary \ref{CorHopf}$)$.

In the last section, as an application of Theorem \ref{MainThm}, we consider some $G$-bivariant analogs of the polynomial ring. In the trivial group case, the polynomial ring $R[\mathbf{X}]$ over $R$ with one variable $\mathbf{X}$ was characterized by the natural bijection
\[ \RAlg(R[\mathbf{X}],S)\overset{\cong}{\rightarrow}S \]
for each $R$-algebra $S$. To any Tambara functor $T$, firstly we associate two types of \lq polynomial' Tambara functors $T[\mathbf{X}]$ and $T[\mathbf{x}]$. Analogously, we obtain natural bijections for each $T$-Tambara functor $S$
\begin{eqnarray*}
\TTamG(T[\mathbf{X}],S)\cong S(G/G), \\
\TTamG(T[\mathbf{x}],S)\cong S(G/e)^G,
\end{eqnarray*}
and thus $T[\mathbf{X}]$ and $T[\mathbf{x}]$ are characterized by these bijections $($Proposition \ref{Prop2}$)$.
Also, $T[\mathbf{X}]$ admits some generalization by using representable Tambara functors (Corollary \ref{CorTX}). On the other hand, $T[\mathbf{x}]$ has a close relationship with the Witt-Burnside ring (Remark \ref{RemCor1WB}).

Finally, as a simultaneous generalization of $T[\mathbf{X}]$ and $T[\mathbf{x}]$, for any subgroup $H\le G$ we construct a functor
\[ \mathit{po\ell}_H\colon\TamG\rightarrow\TamG\ ;\ T\rightarrow T{[}\mathfrak{X}_H{]}, \]
which satisfies the following.
\begin{thmlast}
Let $G$ be a finite group, and $T$ be any Tambara functor on $G$. For any $H\le G$, the Tambara functor $\mathit{po\ell}_H(T)=T[\mathfrak{X}_H]$ satisfies the following.
\begin{itemize}
\item[$(\ast)$] For any $T$-Tambara functor $S$, there exists a natural bijection
\[ \TTamG(T{[}\mathfrak{X}_H{]},S)\cong S(G/H)^{N_G(H)/H}. \]
\end{itemize}
Here, $N_G(H)\le G$ denotes the normalizer of $H$ in $G$.
\end{thmlast}

\bigskip

Throughout this article, we fix a finite group $G$, whose unit element is denoted by $e$.  Abbreviately we denote the trivial subgroup of $G$ by $e$, instead of $\{ e\}$.
$H\le G$ means $H$ is a subgroup of $G$.
$\Gs$ denotes the category of finite $G$-sets and $G$-equivariant maps.
A monoid is always assumed to be unitary and commutative. Similarly a (semi-)ring is assumed to be commutative, and have an additive unit $0$ and a multiplicative unit $1$.
We denote the category of sets by $\Sett$, the category of monoids by $\Mon$, the category of (resp. semi-)rings by $\Ring$ (resp. $\SRing$), and the category of abelian groups by $\Ab$. A monoid homomorphism preserves units, and a (semi-)ring homomorphism preserves $0$ and $1$.

For any category $\mathcal{K}$ and any pair of objects $X$ and $Y$ in $\mathcal{K}$, the set of morphisms from $X$ to $Y$ in $\mathcal{K}$ is denoted by $\mathcal{K}(X,Y)$. For each $X\in\Ob(\mathcal{K})$, the slice category of $\mathcal{K}$ over $X$ is denoted by $\mathcal{K}/X$.

\bigskip

\begin{dfn}\label{DefAddFtr}
An {\it additive contravariant functor} $F$ {\it on} $G$ means a contravariant functor
\[ F\colon\Gs\rightarrow\Sett, \]
which sends coproducts in $\Gs$ to products in $\Sett$.
\end{dfn}

\begin{dfn}\label{DefMackFtr}
A {\it semi-Mackey functor} $M$ {\it on} $G$ is a pair $M=(M^{\ast},M_{\ast})$ of a covariant functor
\[ M_{\ast}\colon\Gs\rightarrow\Sett \]
and an additive contravariant functor
\[ M^{\ast}\colon \Gs\rightarrow\Sett, \]
satisfying $M^{\ast}(X)=M_{\ast}(X)$ for any $X\in\Ob(\Gs)$, and the following Mackey condition:
\begin{itemize}
\item (Mackey condition)

If we are given a pull-back diagram
\[
\xy
(-6,6)*+{X\ppr}="0";
(6,6)*+{Y\ppr}="2";
(-6,-6)*+{X}="4";
(6,-6)*+{Y}="6";
(0,0)*+{\square}="8";
{\ar^{f\ppr} "0";"2"};
{\ar_{x} "0";"4"};
{\ar^{y} "2";"6"};
{\ar_{f} "4";"6"};
\endxy
\]
in $\Gs$, then
\[
\xy
(-12,7)*+{M(X\ppr)}="0";
(12,7)*+{M(Y\ppr)}="2";
(-12,-7)*+{M(X)}="4";
(12,-7)*+{M(Y)}="6";
{\ar_{M^{\ast}(f\ppr)} "2";"0"};
{\ar_{M_{\ast}(x)} "0";"4"};
{\ar^{M_{\ast}(y)} "2";"6"};
{\ar^{M^{\ast}(f)} "6";"4"};
{\ar@{}|\circlearrowright "0";"6"};
\endxy
\]
is commutative.
\end{itemize}
Here we put $M(X)=M^{\ast}(X)=M_{\ast}(X)$ for each $X\in\Ob(\Gs)$. Those $M_{\ast}(f)$ and $M^{\ast}(f)$ for morphisms $f$ in $\Gs$ are called {\it structure morphisms} of $M$. For each $f\in\Gs(X,Y)$, $M^{\ast}(f)$ is called the {\it restriction}, and $M_{\ast}(f)$ is called the {\it transfer} along $f$. If $M$ is a semi-Mackey functor, then $M^{\ast}$ and $M_{\ast}$ become naturally functors to $\Mon$.

For semi-Mackey functors $M$ and $N$, a {\it morphism} from $M$ to $N$ is a family of monoid homomorphisms
\[ \varphi=\{\varphi_X\colon M(X)\rightarrow N(X)\}_{X\in\Ob(\Gs)}, \]
natural with respect to both of the contravariant and the covariant parts. The category of semi-Mackey functors is denoted by $\SMackG$.

If $M(X)$ is an abelian group for each $X\in\Ob(\Gs)$, namely if $M^{\ast}$ and $M_{\ast}$ are functors to $\Ab$, then a semi-Mackey functor $M=(M^{\ast},M_{\ast})$ is called a {\it Mackey functor}.
The full subcategory of Mackey functors in $\SMackG$ is denoted by $\MackG$.
\end{dfn}

\begin{ex}\label{ExMackFtr}
$\ \ $
\begin{enumerate}
\item The Burnside ring functor $\Omega\in\Ob(\MackG)$ (see for example \cite{Bouc}, \cite{O-Y3} or \cite{N_TamMack}).
\item If $Q$ is a (not necessarily finite) $G$-monoid, then the correspondence
\[ \mathcal{P}_Q(X)=\{ G\text{-equivariant maps from}\ X\ \text{to}\ Q\} \]
forms a semi-Mackey functor $\mathcal{P}_Q\in\Ob(\SMackG)$ with structure morphisms defined by
\[ \mathcal{P}_Q^{\ast}(f)\colon\mathcal{P}_Q(Y)\rightarrow\mathcal{P}_Q(X)\ ;\ \beta\mapsto\beta\circ f \]
\[ \ \ (\mathcal{P}_Q)_{\ast}(f)\colon\mathcal{P}_Q(X)\rightarrow\mathcal{P}_Q(Y)\ ;\ \alpha\mapsto(Y\ni y\mapsto \prod_{x\in f^{-1}(y)}\!\!\!\!\!\!\alpha(x)\in Q) \]
for each $f\in\Gs(X,Y)$, where $\prod$ denotes the multiplication of elements in $Q$. This $\mathcal{P}_Q$ is called the {\it fixed point functor} associated to $Q$ (see for example \cite{Tam} or \cite{N_TamMack}). If we denote the category of $G$-monoids by $\GMon$, this construction gives a fully faithful functor
\[ \mathcal{P}\colon\GMon\rightarrow\SMackG\ ;\ Q\mapsto\mathcal{P}_Q. \]
Thus $\GMon$ can be regarded as a full subcategory of $\SMackG$ through $\mathcal{P}$.
\end{enumerate}
\end{ex}

\begin{rem}\label{RemMackTensor}
For each pair of Mackey functors $M$ and $N$, its tensor product $M\underset{\Omega}{\otimes}N$ is defined $($also denoted by $M\widehat{\otimes}N$ in \cite{Bouc}$)$, and $\MackG$ becomes a symmetric monoidal category with this tensor product and the unit $\Omega$.

The category of monoids in $\MackG$ is denoted by $\GreenG$, and a monoid $A$ in $\MackG$ is called a {\it Green functor on} $G$. 

By definition of the tensor product (\cite{Bouc}), for each $X\in\Ob(\Gs)$
\[ (M\underset{\Omega}{\otimes}N)(X)=(\underset{A\overset{p}{\rightarrow}X}{\bigoplus}M(A)\underset{\mathbb{Z}}{\otimes}N(A))/\,\mathcal{I}, \]
where $A\overset{p}{\rightarrow}X$ runs over the objects in $\Gs/X$, and $\mathcal{I}$ is the submodule generated by the elements
\begin{eqnarray*}
&M^{\ast}(a)(s\ppr)\otimes t-s\ppr\otimes N_{\ast}(a)(t),\ \ M_{\ast}(a)(s)\otimes t\ppr-s\otimes N^{\ast}(a)(t\ppr) &\\
&(a\!\in\!\Gs/X((A\overset{p}{\rightarrow}X),(A\ppr\overset{p\ppr}{\rightarrow}X)),s\in M(A),t\in N(A),s\ppr\in M(A\ppr),t\ppr\in N(A\ppr)).&
\end{eqnarray*}
In the component of $A\overset{p}{\rightarrow}X$, the image of $s\otimes t$ in $(M\underset{\Omega}{\otimes}N)(X)$ is denoted by $[s\otimes t]_{(A,p)}$ for each $s\in M(A)$ and $t\in N(A)$.

A priori an element $\omega$ in $(M\underset{\Omega}{\otimes}N)(X)$ is a finite sum of the elements of the above form
\[ \omega=\sum_{1\le i\le n}[s_i\otimes t_i]_{(A_i,p_i)}, \]
however $\omega$ can be written by one such an element. In fact, if we put
\[ A=\coprod_{1\le i\le n}A_i,\ \ p=\underset{1\le i\le n}{\bigcup}p_i,\ \ \iota_i\colon A_i\hookrightarrow A\ \ (\text{inclusion}) \]
and put
\[ s=\sum_{1\le i\le n}(\iota_i)_{\ast}(s_i),\ \ t=\sum_{1\le i\le n}(\iota_i)_{\ast}(t_i), \]
then we have
\[ [s\otimes t]_{(A,p)}=\sum_{1\le i\le n}[s_i\otimes t_i]_{(A_i,p_i)}. \]
\end{rem}

\begin{dfn}
For each $f\in\Gs(X,Y)$ and $p\in\Gs(A,X)$, the {\it canonical exponential diagram} generated by $f$ and $p$ is the commutative diagram
\[
\xy
(-14,6)*+{X}="0";
(-14,-6)*+{Y}="2";
(-1,6)*+{A}="4";
(9,6)*+{}="5";
(18,4.7)*+{X\underset{Y}{\times}\Pi_{f}(A)}="6";
(16,2)*+{}="7";
(16,-6)*+{\Pi_{f}(A)}="8";
(0,0)*+{\mathit{exp}}="10";
{\ar_{f} "0";"2"};
{\ar_{p} "4";"0"};
{\ar_>>>>>{e} "5";"4"};
{\ar^>>>>{f\ppr} "7";"8"};
{\ar^{\pi} "8";"2"};
\endxy
\]
where
\[\Pi_{f}(A)=\Set{(y,\sigma)|%
\begin{array}{l}%
y\in Y, \\
\sigma\colon f^{-1}(y)\rightarrow A \ \, \text{is a map of sets},\\
p\circ \sigma\ \, \text{is identity on}\ f^{-1}(y)%
\end{array}},
\]
\[\pi(y,\sigma)=y,\ \ \ \ \ e(x,(y,\sigma))=\sigma(x), \]
and $f\ppr$ is the pull-back of $f$ by $\pi$.
A diagram in $\Gs$ isomorphic to one of the canonical exponential diagrams is called an {\it exponential diagram}.
For the properties of exponential diagrams, see \cite{Tam}.
\end{dfn}

\begin{dfn}\label{DefTamFtr}
A {\it semi-Tambara functor} $T$ {\it on} $G$ is a triplet $T=(T^{\ast},T_+,T_{\bullet})$ of two covariant functors
\[ T_+\colon\Gs\rightarrow\Sett,\ \ T_{\bullet}\colon\Gs\rightarrow\Sett \]
and one additive contravariant functor
\[ T^{\ast}\colon\Gs\rightarrow\Sett \]
which satisfies the following.
\begin{enumerate}
\item $T^{\alpha}=(T^{\ast},T_+)$ and $T^{\mu}=(T^{\ast},T_{\bullet})$ are objects in $\SMackG$. $T^{\alpha}$ is called the {\it additive part} of $T$, and $T^{\mu}$ is called the {\it multiplicative part} of $T$.
\item (Distributive law)
If we are given an exponential diagram
\[
\xy
(-12,6)*+{X}="0";
(-12,-6)*+{Y}="2";
(0,6)*+{A}="4";
(12,6)*+{Z}="6";
(12,-6)*+{B}="8";
(0,0)*+{exp}="10";
{\ar_{f} "0";"2"};
{\ar_{p} "4";"0"};
{\ar_{\lambda} "6";"4"};
{\ar^{\rho} "6";"8"};
{\ar^{q} "8";"2"};
\endxy
\]
in $\Gs$, then
\[
\xy
(-18,7)*+{T(X)}="0";
(-18,-7)*+{T(Y)}="2";
(0,7)*+{T(A)}="4";
(18,7)*+{T(Z)}="6";
(18,-7)*+{T(B)}="8";
{\ar_{T_{\bullet}(f)} "0";"2"};
{\ar_{T_+(p)} "4";"0"};
{\ar^{T^{\ast}(\lambda)} "4";"6"};
{\ar^{T_{\bullet}(\rho)} "6";"8"};
{\ar^{T_+(q)} "8";"2"};
{\ar@{}|\circlearrowright "0";"8"};
\endxy
\]
is commutative.
\end{enumerate}

If $T=(T^{\ast},T_+,T_{\bullet})$ is a semi-Tambara functor, then $T(X)$ becomes a semi-ring for each $X\in\Ob(\Gs)$, whose additive (resp. multiplicative) monoid structure is induced from that on $T^{\alpha}(X)$ (resp. $T^{\mu}(X)$).
Those $T^{\ast}(f), T_+(f),T_{\bullet}(f)$ for morphisms $f$ in $\Gs$ are called {\it structure morphisms} of $T$. For each $f\in\Gs(X,Y)$,
\begin{itemize}
\item $T^{\ast}(f)\colon T(Y)\rightarrow T(X)$ is a semi-ring homomorphism, called the {\it restriction} along $f$. 
\item $T_+(f)\colon T(X)\rightarrow T(Y)$ is an additive homomorphism, called the {\it additive transfer} along $f$.
\item $T_{\bullet}(f)\colon T(X)\rightarrow T(Y)$ is a multiplicative homomorphism, called the {\it multiplicative transfer} along $f$.
\end{itemize}
$T^{\ast}(f),T_+(f),T_{\bullet}(f)$ are often abbreviated to $f^{\ast},f_+,f_{\bullet}$.

A {\it morphism} of semi-Tambara functors $\varphi\colon T\rightarrow S$ is a family of semi-ring homomorphisms
\[  \varphi=\{\varphi_X\colon T(X)\rightarrow S(X) \}_{X\in\Ob(\Gs)}, \]
natural with respect to all of the contravariant and the covariant parts. We denote the category of semi-Tambara functors by $\STamG$.

If $T(X)$ is a ring for each $X\in\Ob(\Gs)$, then a semi-Tambara functor $T$ is called a {\it Tambara functor}. The full subcategory of Tambara functors in $\STamG$ is denoted by $\TamG$.
\end{dfn}

\begin{rem}\label{RemTamFtr}
A semi-Tambara functor $T$ is a Tambara functor if and only if $T^{\alpha}$ is a Mackey functor. Taking the additive parts and the multiplicative parts, we obtain functors
\begin{eqnarray*}
(-)^{\alpha}\colon\TamG\rightarrow\MackG, \\
(-)^{\mu}\colon\TamG\rightarrow\SMackG.
\end{eqnarray*}
(In fact, $(-)^{\alpha}$ factors through the category $\GreenG$. For this, see \cite{N_TamMack} or \cite{Tam_Manu}.)

If $G$ is trivial, these functors are nothing other than the forgetful functors
\begin{eqnarray*}
&(-)^{\alpha}\colon\Ring\rightarrow\Ab,& \\
&(-)^{\mu}\colon\Ring\rightarrow\Mon,&
\end{eqnarray*}
where, for each ring $R$, its image $R^{\alpha}$ $($resp. $R^{\mu}$$)$ is the underlying additive $($resp. multiplicative$)$ monoid of $R$.
\end{rem}

\begin{ex}\label{ExTamFtr}
The Burnside ring functor $\Omega$ is the initial object in $\TamG$. $\Omega$ can be regarded as the $G$-bivariant analog of $\mathbb{Z}$.
\end{ex}

The following is shown in \cite{Tam_Manu}.
\begin{rem}[ \S 12 in \cite{Tam_Manu}]\label{RemTamManu}
If $T$ and $S$ are Tambara functors, then so is $T\underset{\Omega}{\otimes}S$. Besides, there exist morphisms $\iota_T\in\TamG(T,T\underset{\Omega}{\otimes}S)$ and $\iota_S\in\TamG(S,T\underset{\Omega}{\otimes}S)$, which make $T\underset{\Omega}{\otimes}S$ the coproduct of $T$ and $S$ in $\TamG$. $\Omega$ is the unit for this tensor product.
\end{rem}
\begin{proof}
A simple proof using a functor category will be found in \cite{Tam_Manu}.
For the later use, we briefly introduce an explicit construction of the structure morphisms of $T\underset{\Omega}{\otimes}S$.

Let $f\in\Gs(X,Y)$ be any morphism.
For any $[v\otimes u]_{(B,q)}\in (T\underset{\Omega}{\otimes}S)(Y)$, we define $f^{\ast}([v\otimes u]_{(B,q)})$ by
\[ f^{\ast}([v\otimes u]_{(B,q)})=[T^{\ast}(p_B)(v)\otimes S^{\ast}(p_B)(u)]_{(X\underset{Y}{\times}B,p_X)}, \]
where
\begin{equation}
\label{EqCanPB}
\xy
(-8,6)*+{X\times_YB}="0";
(8,6)*+{B}="2";
(-8,-6)*+{X}="4";
(8,-6)*+{Y}="6";
(0,0)*+{\square}="8";
{\ar^<<<{p_B} "0";"2"};
{\ar_{p_X} "0";"4"};
{\ar^{q} "2";"6"};
{\ar_{f} "4";"6"};
\endxy
\end{equation}
is the canonical pull-back.

For any $[t\otimes s]_{(A,p)}\in (T\underset{\Omega}{\otimes} S)(X)$, we define $f_+([t\otimes s]_{(A,p)})$ and $f_{\bullet}([t\otimes s]_{(A,p)})$ by
\[ f_+([t\otimes s]_{(A,p)})=[t\otimes s]_{(A,f\circ p)}, \]
\[ f_{\bullet}([t\otimes s]_{(A,p)})=[T_{\bullet}(f\ppr)T^{\ast}(e)(t)\otimes S_{\bullet}(f\ppr)S^{\ast}(e)(s)]_{(\Pi_f(A),\pi)}, \]
where
\begin{equation}
\label{EqCanExp}
\xy
(-14,6)*+{X}="0";
(-14,-6)*+{Y}="2";
(-1,6)*+{A}="4";
(9,6)*+{}="5";
(18,4.7)*+{X\underset{Y}{\times}\Pi_{f}(A)}="6";
(16,2)*+{}="7";
(16,-6)*+{\Pi_{f}(A)}="8";
{\ar_{f} "0";"2"};
{\ar_{p} "4";"0"};
{\ar_>>>>>{e} "5";"4"};
{\ar^>>>>{f\ppr} "7";"8"};
{\ar^{\pi} "8";"2"};
{\ar@{}|\circlearrowright "0";"8"};
\endxy
\end{equation}
is the canonical exponential diagram.
With these structure morphisms, $T\underset{\Omega}{\otimes}S$ becomes a Green functor as shown in \cite{Bouc}. Moreover for these (to-be-)structure morphisms, the functoriality of $(T\underset{\Omega}{\otimes}S)_{\bullet}$ follows from $(1.3)$ in \cite{Tam}, the Mackey condition for $(T\underset{\Omega}{\otimes}S)^{\mu}$ follows from $(1.1)$ in \cite{Tam}, the distributive law for $T\underset{\Omega}{\otimes}S$ follows from $(1.2)$ in \cite{Tam}, and thus $T\underset{\Omega}{\otimes}S$ becomes a Tambara functor.
\end{proof}

\begin{dfn}\label{DefT-Tam}
Fix a Tambara functor $T\in\Ob(\TamG)$. 
A $T$-{\it Tambara functor}
is a pair $(S,\sigma)$ of a Tambara functor $S$ and $\sigma\in\TamG(T,S)$. We often represent $(S,\sigma)$ merely by $S$. 

If $(S,\sigma)$ and $(S\ppr,\sigma\ppr)$ are $T$-Tambara functors, then a {\it morphism} $\varphi$ from $(S,\sigma)$ to $(S\ppr,\sigma\ppr)$ is a morphism $\varphi\in\TamG(S,S\ppr)$ satisfying $\sigma\ppr=\varphi\circ\sigma$.
The category of $T$-Tambara functors is denoted by $\TTamG$.
Remark that $\Omega\text{-}\TamG$ is nothing other than $\TamG$.

By Remark \ref{RemTamManu}, for any $T,S\in\Ob(\TamG)$, their tensor product $T\underset{\Omega}{\otimes}S$ can be naturally regarded as a $T$-Tambara functor $(T\underset{\Omega}{\otimes}S,\iota_T)$, or a $S$-Tambara functor $(T\underset{\Omega}{\otimes}S,\iota_S)$.
\end{dfn}

\section{Generalization of the Dress construction}
\label{Section2}
The {\it Dress construction} of a Tambara functor is a process making a Tambara functor $T_Q$ out of a Tambara functor $T$ and a finite $G$-monoid $Q$. This is realized as a functor
\[ \TamG\times\Gmon\rightarrow\TamG \]
where $\Gmon$ is the category of {\it finite} $G$-monoids. Remark here $\Gmon$ is a full subcategory of $\GMon$, and thus can be regarded as a full subcategory of $\SMackG$ through $\mathcal{P}$ (Example \ref{ExMackFtr}).

In this section, we extend this functor to
\[ \TamG\times\SMackG\rightarrow\TamG, \]
through a $G$-bivariant analogical construction of a semi-group ring with a coefficient, by means of the {\it Tambarization functor}.
First, we recall the Dress construction for a Tambara functor:

\begin{dfn}[Theorem 2.9 in \cite{O-Y3}]\label{DefDressTam}
Let $Q$ be a finite $G$-monoid, and $T$ be a Tambara functor on $G$. In \cite{O-Y3}, $T_Q\in\Ob(\TamG)$ is defined by $T_Q(X)=T(X\times Q)$ for each $X\in\Ob(\Gs)$, whose structure morphisms are given by
\begin{eqnarray*}
(T_Q)^{\ast}(f)&=&T^{\ast}(f\times Q)\\
(T_Q)_+(f)&=&T_+(f\times Q)\\
(T_Q)_{\bullet}(f)&=&T_+(\mu_f)\circ T_{\bullet}(f\ppr)\circ T^{\ast}(e)
\end{eqnarray*}
for each morphism $f\in\Gs(X,Y)$.
Here,
\[
\xy
(-24,6)*+{X}="0";
(-24,-6)*+{Y}="2";
(-4,6)*+{X\times Q}="4";
(24,6)*+{X\times_Y\Pi_f(X\times Q)}="6";
(24,-6)*+{\Pi_f(X\times Q)}="8";
(0,0)*+{\mathit{exp}}="10";
{\ar_{f} "0";"2"};
{\ar_{p_X} "4";"0"};
{\ar_<<<<{e} "6";"4"};
{\ar^{f\ppr} "6";"8"};
{\ar^{\pi} "8";"2"};
\endxy
\]
is the canonical exponential diagram with $p_X\colon X\times Q\rightarrow X$ the projection, and $\mu_f\colon\Pi_f(X\times Q)\rightarrow Y\times Q$ is the morphism defined by
\[ \mu_f(y,\sigma)=(y,\prod_{x\in f^{-1}(y)}\! p_Q\circ \sigma(x)), \]
where $p_Q\colon X\times Q\rightarrow Q$ is the projection.

Especially if $T=\Omega$, then $\Omega_Q$ is called the {\it crossed Burnside ring functor}.
\end{dfn}

The crossed Burnside ring functors were generalized by the following Theorems.
\begin{fact}[Theorem 2.15 in \cite{N_TamMack}]\label{FactNTM1}
$(-)^{\mu}\colon\TamG\rightarrow\SMackG$ has a left adjoint functor
\[ \mathcal{T}\colon\SMackG\rightarrow\TamG, \]
which we call the {\it Tambarization functor}.
\end{fact}
\begin{proof}
We briefly review the structure of $\mathcal{T}(M)$. For the entire proof, see \cite{N_TamMack}.

For each $X\in\Ob(\Gs)$, we define $\mathcal{T}(M)$ to be the Grothendieck ring of the category of pairs $(A\overset{p}{\rightarrow}X,m_A)$ of $(A\overset{p}{\rightarrow}X)\in\Ob(\Gs/X)$ and $m_A\in M(A)$.

For each $f\in\Gs(X,Y)$, the structure morphisms induced from $f$ are those determined by
\[ \mathcal{T}(M)^{\ast}(f)(B\overset{q}{\rightarrow}Y,m_B)%
=(X\times_YB\overset{p_X}{\rightarrow}X,M^{\ast}(p_B)(m_B)) \]
\[ \mathcal{T}(M)_+(A\overset{p}{\rightarrow}X,m_A)%
=(A\overset{f\circ p}{\longrightarrow}Y,m_A) \]
\[ \mathcal{T}(M)_{\bullet}(A\overset{p}{\rightarrow}X,m_A)%
=(\Pi_f(A)\overset{\pi}{\rightarrow}Y,M_{\ast}(f\ppr)M^{\ast}(e)(m_A)) \]
for each $(A\overset{p}{\rightarrow}X,m_A)$ and $(B\overset{q}{\rightarrow}Y,m_B)$, where $(\ref{EqCanPB})$ is the canonical pull-back, and $(\ref{EqCanExp})$ is the canonical exponential diagram.
\end{proof}

\begin{fact}[Proposition 3.2 in \cite{N_TamMack}]\label{FactNTM2}
If $Q$ is a finite $G$-monoid, then we have an isomorphism of Tambara functors $\Omega_Q\cong\mathcal{T}(\mathcal{P}_Q)$.
\end{fact}

Thus $\mathcal{T}(M)$, where $M$ can be taken as an arbitrary semi-Mackey functor on $G$, is regarded as a generalization of the crossed Burnside ring functors.

\begin{rem}\label{RemNTM1}
By the adjoint property in Fact \ref{FactNTM1}, $\mathcal{T}(M)$ can be also regarded as a $G$-bivariant analog of the semi-ring. In fact if $G$ is trivial, then $\mathcal{T}$ is nothing other than the functor taking semi-group rings
\[ \mathbb{Z}[-]\colon\Mon\rightarrow\Ring. \]
In this view, from here we denote $\mathcal{T}(M)$ by $\Omega[M]$ instead, for each semi-Mackey functor $M$ on $G$.
\end{rem}

\smallskip

Our aim
is to show the following, which unifies the Tambarization and the Dress construction.
\begin{thm}\label{MainThm}
For any finite group $G$, there is a functor
\[ \mathcal{F}\colon\TamG\times\SMackG\rightarrow\TamG \]
which satisfies the following.
\begin{enumerate}
\item[{\rm (i)}] If $G$ is trivial, then $\mathcal{F}$ agrees with the functor taking semi-group rings with coefficients
\[ \Ring\times\Mon\rightarrow\Ring\ ;\ (R,Q)\mapsto R[Q]. \]
\item[{\rm (ii)}] If $T=\Omega$, we have a natural isomorphism $\mathcal{F}(\Omega,M)\cong\Omega[M]$ for each semi-Mackey functor $M$. 
\item[{\rm (iii)}] If $Q$ is a finite $G$-monoid, then we have a natural isomorphism $\mathcal{F}(T,\mathcal{P}_Q)\cong T_Q$ for each Tambara functor $T$. 
\item[{\rm (iv)}] For each $T$ and $M$, naturally $\mathcal{F}(T,M)$ becomes a $T$-Tambara functor. Moreover if we fix a Tambara functor $T$, then the induced functor
\[ \mathcal{F}(T,-)\colon\SMackG\rightarrow\TTamG \]
is left adjoint to the composition of forgetful functors
\[ \TTamG\rightarrow\TamG\overset{(-)^{\mu}}{\longrightarrow}\SMackG\ ;\ (S,\sigma)\mapsto S^{\mu}. \]
\end{enumerate}
\end{thm}
\begin{proof}
As a consequence of Remark \ref{RemTamManu}, we obtain a functor
\[ \TamG\times\TamG\rightarrow\TamG\ ;\ (T,S)\mapsto T\underset{\Omega}{\otimes}S. \]

Combining this with the Tambarization functor, we define $\mathcal{F}$ by
\[ \mathcal{F}\colon\TamG\times\SMackG\rightarrow \TamG\ ;\ (T,M)\mapsto T\underset{\Omega}{\otimes}\Omega[M]. \]

Then properties {\rm (i)} and {\rm (iv)} follows immediately from the adjoint property of the Tambarization functor, and the universality of the coproduct.
Also, {\rm (ii)} follows from the unit isomorphism $\Omega\underset{\Omega}{\otimes}\Omega[M]\cong\Omega[M]$. Thus it remains to show {\rm (iii)}. 

To show {\rm (iii)}, let $T$ be a Tambara functor, and let $Q$ be a finite $G$-monoid. We construct a natural isomorphism $T_Q\cong T\underset{\Omega}{\otimes}\Omega[\mathcal{P}_Q]$ of Tambara functors. Remark that for each $X\in\Ob(\Gs)$, any element $\omega$ in $(T\underset{\Omega}{\otimes}\Omega[\mathcal{P}_Q])(X)$ can be written in the form of
\begin{equation}
\omega=[s\otimes(R\overset{r}{\rightarrow}A,m_R)]_{(A,p)},
\label{Eq_alpha}
\end{equation}
where $s\in T(A)$, $(A\overset{p}{\rightarrow}X)\in\Ob(\Gs/X)$, $r\in\Gs(R,A)$ and $m_R\in\Gs(R,Q)$.

For each $X\in\Ob(\Gs)$, define
\begin{eqnarray*}
\varphi_X\colon T_Q(X)\rightarrow (T\underset{\Omega}{\otimes}\Omega[\mathcal{P}_Q])(X)\\
\psi_X\colon (T\underset{\Omega}{\otimes}\Omega[\mathcal{P}_Q])(X)\rightarrow T_Q(X)
\end{eqnarray*}
by
\begin{eqnarray*}
&\varphi_X(t)=[t\otimes(X\times Q\overset{\mathrm{id}}{\rightarrow}X\times Q,p_Q)]_{(X\times Q,p_X)}\quad ({}^{\forall}t\in T_Q(X)),&\\
&\psi_X(\omega)=T_+((p\circ r,m_R))T^{\ast}(r)(s)\quad({}^{\forall}\omega\ \text{as in}\ (\ref{Eq_alpha})).&
\end{eqnarray*}
Then we have
\[ \psi_X\circ\varphi_X(t)=T_+((p_X\circ\mathrm{id}_{X\times Q},p_Q))T^{\ast}(\mathrm{id})(t)=t \]
for any $t$, and
\begin{eqnarray*}
\varphi_X\circ\psi_X(\omega)&=&[T_+((p\circ r,m_R))T^{\ast}(r)(s)\otimes(X\times Q\overset{\mathrm{id}}{\rightarrow}X\times Q,p_Q)]_{(X\times Q,p_X)}\\
&=&[s\otimes\Omega[\mathcal{P}_Q]_+(r)\Omega[\mathcal{P}_Q]^{\ast}((p\circ r,m_R))(X\times Q\overset{\mathrm{id}}{\rightarrow}X\times Q,p_Q)]_{(A,p)}\\
&=&[s\otimes(R\overset{r}{\rightarrow}A,m_R)]_{(A,p)}\,=\,\omega
\end{eqnarray*}
for any $\omega$,
namely $\varphi_X$ and $\psi_X$ are mutually inverses.
Thus it remains to show the following:
\begin{claim}
$\varphi=\{\varphi_X\}_{X\in\Ob(\Gs)}$ gives a morphism of Tambara functors
\[ \varphi\colon T_Q\rightarrow T\underset{\Omega}{\otimes}\Omega[\mathcal{P}_Q]. \]
\end{claim}
\begin{proof}
Let $f\in\Gs(X,Y)$ be any morphism. It suffices to show $\varphi$ is compatible with $f^{\ast}$, $f_+$ and $f_{\bullet}$. Denote the projections onto $Q$ by $p_Q\colon X\times Q\rightarrow Q$ and $p_Q\ppr\colon Y\times Q\rightarrow Q$.

\bigskip

\noindent{\rm (a)}\,(compatibility with $f^{\ast}$)

Let $u$ be any element in $T_Q(Y)$. Then we have
\begin{eqnarray*}
f^{\ast}\varphi_Y(u)&=&f^{\ast}([u\otimes(Y\times Q\overset{\mathrm{id}}{\rightarrow}Y\times Q,p_Q\ppr)]_{(Y\times Q,p_Y)})\\
&=&[T^{\ast}(f\times Q)(u)\otimes(X\times Q\overset{\mathrm{id}}{\rightarrow}X\times Q,p_Q)]_{(X\times Q,p_X)}\\
&=&\varphi_X(T^{\ast}(f\times Q)(u))\\
&=&\varphi_XT_Q^{\ast}(f)(u).
\end{eqnarray*}

\smallskip

\noindent{\rm (b)}\,(compatibility with $f_+$)

Let $t$ be any element in $T_Q(X)$. Then we have
\begin{eqnarray*}
f_+\varphi_X(t)&=&f_+([t\otimes(X\times Q\overset{\mathrm{id}}{\rightarrow}X\times Q)]_{(X\times Q,p_X)})\\
&=&[t\otimes(X\times Q\overset{\mathrm{id}}{\rightarrow}X\times Q)]_{(X\times Q,f\circ p_X)}\\
&=&[t\otimes(\Omega[\mathcal{P}_Q]^{\ast}(f\times Q))(Y\times Q\overset{\mathrm{id}}{\rightarrow}Y\times Q,p_Q\ppr)]_{(X\times Q,f\circ p_X)}\\
&=&[T_+(f\times Q)(t)\otimes(Y\times Q\overset{\mathrm{id}}{\rightarrow}Y\times Q,p_Q\ppr)]_{(X\times Q,f\circ p_X)}\\
&=&\varphi_Y(T_+(f\times Q)(t))\\
&=&\varphi_YT_{Q+}(f)(t).
\end{eqnarray*}

\smallskip

\noindent{\rm (c)}\,(compatibility with $f_{\bullet}$)

We use the notation in Definition \ref{DefDressTam}.
For any element $t$ in $T_Q(X)$, we have
\begin{eqnarray*}
\varphi_YT_{Q\bullet}(t)&=&\varphi_Y(T_+(\mu_f)T_{\bullet}(f\ppr)T^{\ast}(e)(t))\\
&=&[T_+(\mu_f)T_{\bullet}(f\ppr)T^{\ast}(e)(t)\otimes(Y\times Q\overset{\mathrm{id}}{\rightarrow}Y\times Q,p_Q\ppr)]_{(Y\times Q,p_Y)}\\
&=&[T_{\bullet}(f\ppr)T^{\ast}(e)(t)\otimes\Omega[\mathcal{P}_Q]^{\ast}(\mu_f)(Y\times Q\overset{\mathrm{id}}{\rightarrow}Y\times Q,p_Q\ppr)]_{(\Pi_f(X\times Q),\pi)}.
\end{eqnarray*}
On the other hand, we have
\begin{eqnarray*}
&f_{\bullet}&\!\!\!\!\!\!\varphi_X(t)= f_{\bullet}([t\otimes(X\times Q\overset{\mathrm{id}}{\rightarrow}X\times Q,p_Q)]_{(X\times Q,p_X)})\\
&=&\!\!\! [T_{\bullet}(f\ppr)T^{\ast}(e)(t)\otimes\Omega[\mathcal{P}_Q]_{\bullet}(f\ppr)\Omega[\mathcal{P}_Q]^{\ast}(e)(X\times Q\overset{\mathrm{id}}{\rightarrow}X\times Q,p_Q)]_{(\Pi_f(X\times Q),\pi)}.
\end{eqnarray*}
Thus it suffices to show
\begin{eqnarray}\label{Eq_Conclusion}
\Omega[\mathcal{P}_Q]^{\ast}(\mu_f)(Y\times Q\!\!\!\!&\overset{\mathrm{id}}{\rightarrow}&\!\!\!\! Y\times Q,p_Q\ppr)\\
&=&\!\!\!\Omega[\mathcal{P}_Q]_{\bullet}(f\ppr)\Omega[\mathcal{P}_Q]^{\ast}(e)(X\times Q\overset{\mathrm{id}}{\rightarrow}X\times Q,p_Q).\nonumber
\end{eqnarray}

The left hand side of $(\ref{Eq_Conclusion})$ is equal to
\begin{equation}
(\Pi_f(X\times Q)\overset{\mathrm{id}}{\rightarrow}\Pi_f(X\times Q),p_Q\ppr\circ\mu_f).
\label{Eq_mult1}
\end{equation}
Remark that $p_Q\ppr\circ\mu_f\colon\Pi_f(X\times Q)\rightarrow Q$ is the morphism which satisfies
\[ \ p_Q\ppr\circ\mu_f(y,\sigma)=p_Q\ppr(y,\!\!\!\! \prod_{x\in f^{-1}(y)}\!\!\! p_Q\circ\sigma(x))=\!\!\!\prod_{x\in f^{-1}(y)}\!\!\! p_Q\circ\sigma(x) \]
for any $(y,\sigma)\in\Pi_f(X\times Q)$.

On the other hand, since there is an exponential diagram
\[
\xy
(-36,7)*+{X\underset{Y}{\times}\Pi_f(X\times Q)}="0";
(-36,-7)*+{\Pi_f(X\times Q)}="2";
(0,7)*+{X\underset{Y}{\times}\Pi_f(X\times Q)}="4";
(36,7)*+{X\underset{Y}{\times}\Pi_f(X\times Q)}="6";
(36,-7)*+{\Pi_f(X\times Q)}="8";
(0,0)*+{\mathit{exp}}="10";
(46,-8)*+{,}="21";
{\ar_{f\ppr} "0";"2"};
{\ar_{\mathrm{id}} "4";"0"};
{\ar_>>>>>{\mathrm{id}} "6";"4"};
{\ar^{f\ppr} "6";"8"};
{\ar^{\mathrm{id}} "8";"2"};
\endxy
\]
the right hand side of $(\ref{Eq_Conclusion})$ is equal to
\[ (\Pi_f(X\times Q)\overset{\mathrm{id}}{\rightarrow}\Pi_f(X\times Q),(\mathcal{P}_Q)_{\ast}(f\ppr)(p_Q\circ e)), \]
where $(\mathcal{P}_Q)_{\ast}(f\ppr)(p_Q\circ e)\colon\Pi_f(X\times Q)\rightarrow Q$ is the morphism satisfying
\begin{equation}
((\mathcal{P}_Q)_{\ast}(f\ppr)(p_Q\circ e))(y,\sigma)=\!\!\! \prod_{(x,(y,\sigma))\in f^{\prime -1}(y,\sigma)}\!\!\! p_Q\circ e(x,(y,\sigma))
\label{Eq_mult2}
\end{equation}
for each $(y,\sigma)\in\Pi_f(X\times Q)$.

Since $f^{\prime -1}(y,\sigma)=\{(x,(y,\sigma))\mid x\in f^{-1}(y) \}$, we have
\begin{eqnarray}
\label{Eq_mult3}
\prod_{(x,(y,\sigma))\in f^{\prime -1}(y,\sigma)}p_Q\circ e(x,(y,\sigma))%
&=&\!\!\! \prod_{(x,(y,\sigma))\in f^{\prime -1}(y,\sigma)}\!\!\! p_Q\circ\sigma(x)\\\nonumber
&=&\!\!\! \prod_{x\in f^{-1}(y)}\!\!\! p_Q\circ \sigma(x)\\\nonumber
\end{eqnarray}
By $(\ref{Eq_mult1})$, $(\ref{Eq_mult2})$ and $(\ref{Eq_mult3})$, we obtain
\[ p_Q\ppr\circ\mu_f=(\mathcal{P}_Q)_{\ast}(f\ppr)(p_Q\circ e), \]
and the equality of $(\ref{Eq_Conclusion})$ follows.
\end{proof}
\end{proof}

\begin{rem}\label{RemMainThm}
By {\rm (i)} and {\rm (iv)} in Theorem \ref{MainThm}, $\mathcal{F}$ can be regarded as a $G$-bivariant analog of the functor taking semi-group rings with coefficients. In this view, from here we denote $\mathcal{F}(T,M)$ by $T[M]$ instead.
\end{rem}

As in the trivial group case, $T[M]$ is equipped with a natural Hopf structure if $M$ is a Mackey functor.
To state this, first we remark the following.
\begin{rem}\label{RemHopf}
For each pair $M,N\in\Ob(\MackG)$, the coproduct of $T[M]$ and $T[N]$ in $\TamG$ is given by $(T[M\oplus N],i_M,i_N)$, where $i_M$ and $i_N$ are the morphisms
\begin{eqnarray*}
&i_M\,\colon\, T[M]\hookrightarrow T[M\oplus N]&\\
&i_N\,\colon\, T[N]\hookrightarrow T[M\oplus N]&
\end{eqnarray*}
induced from the inclusions of semi-Mackey functors $M\hookrightarrow M\oplus N,\, N\hookrightarrow M\oplus N$.
\end{rem}
\begin{proof}
This immediately follows from the adjointness shown in Theorem \ref{MainThm}.
\end{proof}

\begin{cor}\label{CorHopf}
For any $M\in\Ob(\MackG)$, there exist morphisms of $T$-Tambara functors
\begin{eqnarray*}
\Delta_M&\colon&T[M]\rightarrow T[M]\underset{T}{\otimes}T[M]\\
\varepsilon_M&\colon&T[M]\rightarrow T\\
\eta_M&\colon&T[M]\rightarrow T[M]
\end{eqnarray*}
satisfying
\[ (\Delta_M\underset{T}{\otimes}\mathrm{id})\circ T=(\mathrm{id}\underset{T}{\otimes}\Delta_M)\circ T,\ \ (\varepsilon_M\underset{T}{\otimes}\mathrm{id})\circ\Delta_M=\mathrm{id}=(\mathrm{id}\underset{T}{\otimes}\varepsilon_M)\circ\Delta_M \]
and
\[ \nabla\circ(\eta_M\underset{T}{\otimes}\mathrm{id})\circ\Delta_M=\varepsilon_M\circ\iota_T=\nabla\circ(\mathrm{id}\underset{T}{\otimes}\eta_M)\circ\Delta_M, \]
where $\nabla$ is the multiplication morphism $($i.e. the morphism inducing $\mathrm{id}_{T[M]}$ on each components of $T[M]\underset{T}{\otimes}T[M]$.$)$
\[
\xy
(-24,6)*+{T[M]}="0";
(24,6)*+{T[M]\otimes_TT[M]}="2";
(-24,-6)*+{T[M]\otimes_TT[M]}="4";
(24,-6)*+{T[M]\otimes_TT[M]\otimes_TT[M]}="6";
{\ar^{\Delta_M} "0";"2"};
{\ar_{\Delta_M} "0";"4"};
{\ar^{\Delta_M\underset{T}{\otimes}\mathrm{id}} "2";"6"};
{\ar_<<<<<<{\mathrm{id}\underset{T}{\otimes}\Delta_M} "4";"6"};
{\ar@{}|\circlearrowright "0";"6"};
\endxy
\]
$
\xy
(-20,6)*+{T[M]}="0";
(20,6)*+{T\otimes_TT[M]}="2";
(4,-4)*+{}="3";
(-20,-6)*+{T[M]\otimes_TT}="4";
(-4,4)*+{}="5";
(20,-6)*+{T[M]\otimes_TT[M]}="6";
{\ar^{\cong} "0";"2"};
{\ar_{\cong} "0";"4"};
{\ar_{\varepsilon_M\underset{T}{\otimes}\mathrm{id}} "6";"2"};
{\ar^{\mathrm{id}\underset{T}{\otimes}\varepsilon_M} "6";"4"};
{\ar^{\Delta_M} "0";"6"};
{\ar@{}|\circlearrowright "2";"3"};
{\ar@{}|\circlearrowright "4";"5"};
\endxy
$
$
\xy
(-16,10)*+{T[M]\otimes_TT[M]}="0";
(16,10)*+{T[M]\otimes_TT[M]}="2";
(-20,0)*+{T[M]}="4";
(0,9)*+{}="5";
(0,0)*+{T}="6";
(0,-9)*+{}="7";
(20,0)*+{T[M]}="8";
(-16,-10)*+{T[M]\otimes_TT[M]}="10";
(16,-10)*+{T[M]\otimes_TT[M]}="12";
{\ar^{\Delta_M} "4";"0"};
{\ar^{\eta_M\underset{T}{\otimes}\mathrm{id}} "0";"2"};
{\ar^{\nabla} "2";"8"};
{\ar^{\varepsilon_M} "4";"6"};
{\ar^{\iota_T} "6";"8"};
{\ar_{\Delta_M} "4";"10"};
{\ar_{\mathrm{id}\underset{T}{\otimes}\eta_M} "10";"12"};
{\ar_{\nabla} "12";"8"};
{\ar@{}|\circlearrowright "6";"5"};
{\ar@{}|\circlearrowright "6";"7"};
\endxy
$
\end{cor}
\begin{proof}
We define $\Delta_M, \varepsilon_M, \eta_M$ to be those morphisms corresponding to 
\[ M\overset{\text{diag.}}{\longrightarrow}M\oplus M,\ \ M\overset{0}{\longrightarrow}0,\ \ M\overset{(-)^{-1}}{\longrightarrow}M \]
which are the diagonal morphism for $M$, the zero morphism, and the inverse, respectively (defined in an obvious way). Then the required compatibility conditions immediately follows from the functoriality of $\mathcal{F}(T,-)$.
\end{proof}

\section{Some remarks on Hopf Tambara functors}

Concerning Corollary \ref{CorHopf}, we investigate a bit more about Hopf $T$-Tambara functors, restricting ourselves to the case $T=\Omega$.
In this case, for any $M,N\in\Ob(\SMackG)$, the Tambara functors $\Omega[M]\underset{\Omega}{\otimes}\Omega[N]$ in Remark \ref{RemTamManu} and $\Omega[M\oplus N]$ become isomorphic thorough the morphism 
\[ c\,\colon\, \Omega{[}M{]}\underset{\Omega}{\otimes}\Omega{[}N{]}\overset{\cong}{\longrightarrow}\Omega{[}M\oplus N{]}, \]
which arising from the inclusions $\Omega[M]\hookrightarrow \Omega[M\oplus N]$ and $\Omega[N]\hookrightarrow \Omega[M\oplus N]$. More explicitly, for each $X\in\Ob(\Gs)$,
\[ c_X\colon\ (\Omega[M]\underset{\Omega}{\otimes}\Omega[N])(X)\longrightarrow(\Omega[M\oplus N])(X) \]
satisfies
\begin{eqnarray*}
\lefteqn{c_X(\, {[}(A\overset{p}{\rightarrow}X,m_A)\otimes (B\overset{q}{\rightarrow}X,m_B){]}_{(X,\id_X)}\, )}\hspace{3cm}\\
&=&(\, A\underset{X}{\times}B\overset{r}{\rightarrow}X,\, (M^{\ast}(\pi_A)(m_A),N^{\ast}(\pi_B)(m_B))\, )
\end{eqnarray*}
for any $(A\overset{p}{\rightarrow}X,m_A)\in\Omega[M](X)$ and $(B\overset{q}{\rightarrow}X,m_B)\in\Omega[N](X)$, where
\[
\xy
(-8,6)*+{A\times_XB}="0";
(8,6)*+{B}="2";
(-8,-6)*+{A}="4";
(8,-6)*+{X}="6";
(0,0)*+{\square}="7";
{\ar^>>>>{\pi_B} "0";"2"};
{\ar_{\pi_A} "0";"4"};
{\ar_{p} "4";"6"};
{\ar^{q} "2";"6"};
\endxy
\]
is the pullback in $\Gs$, and $r=p\circ\pi_A=q\circ\pi_B$.

\begin{dfn}\label{DefHopf}
A Tambara functor $\HH$ is a {\it Hopf Tambara functor} if it is equipped with morphisms in $\TamG$

\[
\Delta\colon\HH\rightarrow \HH\underset{\Omega}{\otimes}\HH,\quad %
\varepsilon\colon\HH\rightarrow T,\quad%
\eta\colon\HH\rightarrow \HH,
\]
satisfying
\begin{eqnarray*}
&(\Delta\underset{\Omega}{\otimes}\mathrm{id})\circ \Delta=(\mathrm{id}\underset{\Omega}{\otimes}\Delta)\circ \Delta,&\\
&(\varepsilon\underset{\Omega}{\otimes}\mathrm{id})\circ\Delta=\mathrm{id}=(\mathrm{id}\underset{\Omega}{\otimes}\varepsilon)\circ\Delta,&\\
&\nabla\circ(\eta\underset{\Omega}{\otimes}\mathrm{id})\circ\Delta=\varepsilon\circ\iota=\nabla\circ(\mathrm{id}\underset{\Omega}{\otimes}\eta)\circ\Delta.&
\end{eqnarray*}

%


If $(\HH,\Delta,\varepsilon,\eta)$ and $(\HH\ppr,\Delta\ppr,\varepsilon\ppr,\eta\ppr)$ are Hopf Tambara functors, then a morphism $\varphi\in\TamG(\HH,\HH\ppr)$ is a {\it morphism of Hopf Tambara functors} if it satisfies
\[
\Delta\ppr\circ \varphi=(\varphi\underset{\Omega}{\otimes}f)\circ\Delta,\quad%
\varepsilon\ppr\circ \varphi=\varepsilon,\quad%
\eta\ppr\circ \varphi=\varphi\circ\eta.
\]

We denote the category of Hopf Tambara functors by $\HopfG$.
\end{dfn}

\begin{rem}
As Corollary \ref{CorHopf} suggests, Tambarization gives a functor
\[ \Omega[-]\colon\MackG\rightarrow\HopfG. \]
\end{rem}

\medskip

As the group-like elements form a multiplicative subgroup in any Hopf algebra over a commutative ring, we can obtain a Mackey subfunctor of a Hopf Tambara functor.
\begin{dfn}\label{DefGroupLike}
Let $\HH$ be a Hopf Tambara functor. For each $X\in\Ob(\Gs)$, we call an element $m\in \HH(X)$ {\it group-like} if it satisfies
\begin{eqnarray}
\label{Eqi}&\Delta_X(m)=[m\otimes m]_{(X,\id_X)}&\quad\text{in}\ \ (\HH\underset{\Omega}{\otimes}\HH)(X),\\
\label{Eqii}&\varepsilon_X(m)=1&\quad\text{in}\ \ \Omega(X).
\end{eqnarray}
\end{dfn}

\begin{prop}\label{PropGroupLike}
Let $\HH$ be a Hopf Tambara functor.
If we define $\M_{\HH}$ by
\[ \M_{\HH}(X)=\{ m\in\HH(X)\mid m\ \text{is group-like} \} \]
for each $X\in\Ob(\Gs)$, then $\M_{\HH}$ becomes a Mackey subfunctor of $\HH^{\mu}$.
\end{prop}
\begin{proof}
First we show $\M_{\HH}(X)$ is a subgroup of $\HH^{\mu}(X)$ for each $X\in\Ob(\Gs)$. Take any $m_1,m_2\in \M_{\HH}(X)$.
From
\begin{eqnarray*}
\Delta_X(m_1m_2)&=&\Delta_X(m_1)\Delta_X(m_2)\\
&=&[m_1\otimes m_1]_{(X,\id_X)}\cdot [m_2\otimes m_2]_{(X,\id_X)}\\
&=&[m_1m_2\otimes m_1m_2]_{(X,\id_X)}
\end{eqnarray*}
and
\[ \varepsilon_X(m_1m_2)=\varepsilon_X(m_1)\varepsilon_X(m_2)=1, \]
we obtain $m_1m_2\in\M_{\HH}(X)$.
Besides, by
\begin{eqnarray*}
1&=&\iota_X\varepsilon_X(m_1)\ =\ \nabla_X(\eta\underset{\Omega}{\otimes}\id)_X\Delta_X(m_1)\\
&=&\nabla_X(\eta\underset{\Omega}{\otimes}\id)_X([m_1\otimes m_1]_{(X,\id_X)})\\
&=&\nabla_X([\eta_X(m_1)\otimes m_1]_{(X,\id_X)})\\
&=&\eta_X(m_1)\cdot m_1,
\end{eqnarray*}
we have $\eta_X(m_1)=m_1^{-1}$. Since we have
\[ 1=\varepsilon_X(1)=\varepsilon(\eta_X(m_1)\cdot m_1)=\varepsilon_X(\eta_X(m_1)), \]
it follows $m_1^{-1}=\eta_X(m_1)\in\M_{\HH}(X)$.

It remains to show that $\M_{\HH}$ satisfies
\begin{eqnarray*}
\HH_{\bullet}(f)(\M_{\HH}(X))\subseteq \M_{\HH}(Y)\\
\HH^{\ast}(f)(\M_{\HH}(Y))\subseteq \M_{\HH}(X)
\end{eqnarray*}
for any $f\in\Gs(X,Y)$.
Let $m\in\M_{\HH}(X)$ be any element. We have
\[ \varepsilon_Yf_{\bullet}(m)=f_{\bullet}\varepsilon_X(m)=1. \]
Since
\[
\xy
(-14,6)*+{X}="0";
(-14,-6)*+{Y}="2";
(-1,6)*+{X}="4";
(14,6)*+{X}="6";
(14,-6)*+{Y}="8";
(0,0)*+{\mathit{exp}}="10";
{\ar_{f} "0";"2"};
{\ar_{\id_X} "4";"0"};
{\ar_>>>>>{\id_X} "6";"4"};
{\ar^>>>>{f} "6";"8"};
{\ar^{\id_Y} "8";"2"};
\endxy
\]
is an exponential diagram, we have
\begin{eqnarray*}
\Delta_Yf_{\bullet}(m)=f_{\bullet}\Delta_X(m)=f_{\bullet}([m\otimes m]_{(X,\id_X)})=[f_{\bullet}(m)\otimes f_{\bullet}(m)]_{(Y,\id_Y)},
\end{eqnarray*}
and thus $f_{\bullet}(m)\in\M_{\HH}(Y)$.

Let $n\in\M_{\HH}(Y)$ be any element. We have
\[ \Delta_Xf^{\ast}(n)=f^{\ast}\Delta_Y(n)=f^{\ast}([n\otimes n]_{(Y,\id_Y)})=[f^{\ast}(n)\otimes f^{\ast}(n)]_{(X,\id_X)}, \]
\[ \varepsilon_Xf^{\ast}(n)=f^{\ast}\varepsilon_Y(n)=1, \]
and thus $f^{\ast}(n)\in\M_{\HH}(X)$.
\end{proof}

\begin{cor}
$\ \ $
\begin{enumerate}
\item There is a functor
\[ \M\colon\HopfG\rightarrow\MackG\ ;\ \HH\mapsto\M_{\HH} \]
and a natural transformation
\[
\xy
(-16,6)*+{\HopfG}="0";
(16,6)*+{\TamG}="2";
(-16,-6)*+{\MackG}="4";
(16,-6)*+{\SMackG}="6";
(25,-7)*+{.}="7";
(-3,0)*+{}="1";
(3,0)*+{}="3";
{\ar^{\text{forgetful}} "0";"2"};
{\ar_{\M} "0";"4"};
{\ar^{(\,)^{\mu}} "2";"6"};
{\ar@{^(->} "4";"6"};
{\ar@{=>} "1";"3"};
\endxy
\]
\item For any Mackey functor $M$, there is a natural inclusion of Mackey functors
\[ M\hookrightarrow\M_{\Omega[M]}. \]
\item For any Hopf Tambara functor $\HH$, there is a natural morphism of Tambara functors
\[ \Omega[\M_{\HH}]\rightarrow\HH. \]
\end{enumerate}
\end{cor}
\begin{proof}
$\ \ $
\begin{enumerate}
\item This follows immediately from Proposition \ref{PropGroupLike}. The natural transformation is given by the inclusion $\M_{\HH}\hookrightarrow\HH^{\mu}$ for each $\HH\in\Ob(\HopfG)$.

\item Remark that, by the adjointness between $\Omega[-]$ and $(\ )^{\mu}$ in Fact \ref{FactNTM1}, we have a unit morphism
\[ u\colon M\rightarrow(\Omega[M])^{\mu}. \]
For each $X\in\Ob(\Gs)$, this is given by
\[ M(X)\rightarrow(\Omega[M])^{\mu}(X)\ \ ;\ \ m\mapsto (X\overset{\id_X}{\rightarrow}X,m)\quad({}^{\forall}m\in M(X)), \]
which is inclusive. 

It remains to show $u_X(m)\in\Omega[M](X)$ is group-like for any $m\in M(X)$.
Remark that $\Delta_M\colon \Omega[M]\rightarrow\Omega[M\oplus M]$ maps
\[ (X\overset{\id_X}{\rightarrow}X,m)\in \Omega[M](X) \]
to
\[ (X\overset{\id_X}{\rightarrow}X,(m,m))\in \Omega[M\oplus M](X). \]
Since $(X\overset{\id_X}{\rightarrow}X,(m,m))$ corresponds to
\[ {[}(X\overset{\id_X}{\rightarrow}X,m)\otimes (X\overset{\id_X}{\rightarrow}X,m){]}_{(X,\id_X)}\in(\Omega[M]\underset{\Omega}{\otimes}\Omega[M])(X)\]
under the isomorphism $c\colon \Omega[M]\underset{\Omega}{\otimes}\Omega[N]\overset{\cong}{\longrightarrow}\Omega[M\oplus N]$, this shows $u_X(m)$ satisfies condition $(\ref{Eqi})$ in Definition \ref{DefGroupLike}. Besides, since $\varepsilon_M\colon \Omega{[}M{]}\rightarrow\Omega$ 
maps $u_X(m)$ to $(X\overset{\id_X}{\rightarrow}X)=1$, condition $(\ref{Eqii})$ in Definition \ref{DefGroupLike} is satisfied. Thus $u_X(m)$ is group-like.

\item By {\rm (1)} and the adjointness between $\Omega[-]\colon\SMackG\rightarrow\TamG$ and $(\ )^{\mu}\colon\TamG\rightarrow\SMackG$,
we obtain a morphism of Tambara functors $\Omega[\M_{\HH}]\rightarrow\HH$.
\end{enumerate}
\end{proof}

\section{Polynomial Tambara functors}
In this section, we consider $G$-bivariant analogs of the polynomial ring, by using Theorem \ref{MainThm}. Remark that, in the trivial group case, the polynomial ring satisfies the following properties.
\begin{rem}\label{RemPoly}
Let $R$ be a ring, and let $R[\mathbf{X}]$ be the polynomial ring over $R$ with one variable. Then we have the following.
\begin{enumerate}
\item (Existence of the indeterminate element)
For any $R$-algebra $S$, we have a natural bijection
\[ \RAlg(R[\mathbf{X}],S)\overset{\cong}{\longrightarrow}S\ ;\ \varphi\mapsto\varphi(\mathbf{X}), \]
where $\RAlg$ denotes the category of $R$-algebras.

\item (Structural isomorphism)
We have a natural isomorphism of rings
\begin{equation}
R[\mathbf{X}]\cong R\underset{\mathbb{Z}}{\otimes}\mathbb{Z}[\mathbf{X}]\cong R\underset{\mathbb{Z}}{\otimes}\mathbb{Z}[\mathbb{N}]
\label{EqPoly2}
\end{equation}
\end{enumerate}
\end{rem}

\smallskip

\subsection{First definitions of polynomial Tambara functors}

$\ \ $

\smallskip

We first propose two types of \lq polynomial' Tambara functors, satisfying analogous properties to those in Remark \ref{RemPoly}. These Tambara functors will be generalized later (Corollary \ref{CorTX}, Theorem \ref{ThmSimulPoly}).
\begin{prop}\label{Prop2}
Let $G$ be a finite group.
\begin{enumerate}
\item There exists a functor
\[ \mathit{po\ell}_{\mathbf{X}}\colon\TamG\rightarrow\TamG\ ;\ T\mapsto T[\mathbf{X}], \]
which admits a natural bijection
\[ \TTamG(T[\mathbf{X}],S)\cong S(G/G) \]
for each $T\in\Ob(\TamG)$ and $S\in\Ob(\TTamG)$.
\item There exists a functor
\[ \mathit{po\ell}_{\mathbf{x}}\colon\TamG\rightarrow\TamG\ ;\ T\mapsto T[\mathbf{x}], \]
which admits a natural bijection
\[ \TTamG(T[\mathbf{x}],S)\cong S(G/e)^G \]
for each $T\in\Ob(\TamG)$ and $S\in\Ob(\TTamG)$.
\end{enumerate}
Moreover if $G$ is trivial, each of these agrees with the functor taking the polynomial ring\ \ $\mathit{po\ell}\colon\Ring\rightarrow\Ring\ ;\ R\mapsto R[\mathbf{X}]$.
\end{prop}

\begin{proof}
If we follow the analogy of $(\ref{EqPoly2})$, we can expect that each of the desired functors is of the form
\[ \mathcal{F}(-,M)\colon\TamG\rightarrow\TamG \]
for some semi-Mackey functor $M$, which can be regarded as a \lq $G$-bivariant analog of $\mathbb{N}$'.

\medskip

To show {\rm (1)}, we first remark the following.
\begin{rem}\label{RemSemiBurn}
For any $X\in\Ob(\Gs)$, the set of isomorphism classes $\mathit{c\ell}(\Gs/X)$ of the category $\Gs/X$ forms a semi-ring.
If we define $\mathfrak{A}$ by $\mathfrak{A}(X)=\mathit{c\ell}(\Gs/X)$, then $\mathfrak{A}$ becomes a semi-Tambara functor on $G$, called the {\it Burnside semi-ring functor}, with appropriately defined structure morphisms.

If we denote the isomorphism class of $(G/G\overset{\mathrm{id}}{\longrightarrow}G/G)$ in $\mathfrak{A}(G/G)$ by $\mathbf{X}$, then we have a natural isomorphism 
\[ \SMackG(\mathfrak{A}^{\alpha},M)\cong M(G/G)\ ;\ \varphi\mapsto\varphi_{G/G}(\mathbf{X}) \]
for any $M\in\SMackG$.
\end{rem}
As a corollary of Theorem \ref{MainThm} and Remark \ref{RemSemiBurn}, we obtain:
\begin{cor}\label{Cor2}
Let $T$ be a Tambara functor on $G$. For any $T$-Tambara functor $S$, we have an isomorphism
\[ \TTamG(T[\mathfrak{A}^{\alpha}],S)\cong S^{\mu}(G/G) \]
which is natural in $S$.
\end{cor}
Thus if we denote $T[\mathfrak{A}^{\alpha}]$ abbreviately by $T[\mathbf{X}]$, 
then $T[\mathbf{X}]$ satisfies the desired property in {\rm (1)}. $\mathit{po\ell}_{\mathbf{X}}$ is given by $\mathit{po\ell}_{\mathbf{X}}=\mathcal{F}(-,\mathfrak{A}^{\alpha})$.

If $G$ is trivial, (and thus $T$ is identified with the ring $R=T(G/G)$,) then $T[\mathbf{X}]$ is naturally isomorphic to the polynomial ring $R[\mathbf{X}]$ over $R$, with the indeterminate element $\mathbf{X}$.

\medskip

To show {\rm (2)}, we remark the following.
\begin{rem}[Claim 3.8 in \cite{N_TamMack}]\label{RemMonAdj}
If $M$ is a semi-Mackey functor on $G$, then $M(G/e)$ carries a natural $G$-monoid structure. The functor taking its $G$-fixed part
\[ \ev\colon\SMackG\rightarrow\Mon\ ;\ M\mapsto M(G/e)^G\]
admits a left adjoint functor
\[ \mathcal{L}^G\colon\Mon\rightarrow\SMackG \ ;\ Q\mapsto\mathcal{L}^G_Q. \]
\end{rem}
Combining this with Theorem \ref{MainThm}, we obtain:
\begin{cor}
Let $T$ be a Tambara functor on $G$. For any monoid $Q$ and any $T$-Tambara functor $S$, we have an isomorphism
\[ \TTamG(T[\mathcal{L}^G_Q], S)\cong\Mon(Q, S^{\mu}(G/e)^G) \]
which is natural in $Q$ and $S$.
\end{cor}
Especially when $Q=\mathbb{N}$, then we obtain a natural bijection
\[ \TTamG(T[\mathcal{L}^G_{\mathbb{N}}],S)\cong S^{\mu}(G/e)^G. \]
Thus if we denote 
$\mathcal{L}^G_{\mathbb{N}}$ by $\mathbf{x}$, then $T[\mathbf{x}]$ satisfies the desired property in {\rm (2)}. $\mathit{po\ell}_{\mathbf{x}}$ is given by $\mathit{po\ell}_{\mathbf{x}}=\mathcal{F}(-,\mathcal{L}^G_{\mathbb{N}})$.

If $G$ is trivial, (and thus $T$ is identified with the ring $R=T(G/e)$,) then $T[\mathbf{x}]$ is naturally isomorphic to the polynomial ring over $R$ with one variable.

\end{proof}

\begin{rem}\label{RemCor1WB}
We remark also that $T[\mathbf{x}]$ is closely related to the Witt-Burnside ring. In fact, we have a natural isomorphism of commutative rings
\[ \Omega[\mathbf{x}](G/G)\cong\mathbb{W}_G(\mathbb{Z}[\mathbf{X}]), \]
where the right hand side is the Witt-Burnside ring of the polynomial ring $\mathbb{Z}[\mathbf{X}]$ over $G$. (Theorem 3.9 in \cite{N_TamMack}, Theorem 1.7 in \cite{Elliott}. See also \cite{Brun}.)
\end{rem}

\smallskip

\subsection{Representability of $\Omega [\mathbf{X}]$}

$\ \ $

\smallskip

We show $\Omega [\mathbf{X}]$ is a \lq representable' functor in some sense. To state this precisely, we recall the following.

\begin{fact}(\cite{Tam})
Let $G$ be a finite group.
\begin{enumerate}
\item There exists a category $\U_G$ with finite products, satisfying the following properties.
\begin{enumerate}
\item[{\rm (i)}] $\Ob(\U_G)=\Ob(\Gs)$.
\item[{\rm (ii)}] For any $X,Y\in\Ob(\U_G)$, the set of morphisms $\U_G(X,Y)$ becomes a semi-ring. 
\item[{\rm (iii)}] There is a categorical equivalence $\mu_G\colon\Add(\U_G,\Sett)\overset{\simeq}{\longrightarrow}\STamG$. 
\end{enumerate}
\item There exists a category $\V_G$ with finite products, satisfying the following.
\begin{enumerate}
\item[{\rm (i)}] $\Ob(\U_G)=\Ob(\V_G)$.
\item[{\rm (ii)}] For any $X,Y\in\Ob(\V_G)$, the set of morphisms $\V_G(X,Y)$ is isomorphic to the ring-completion $K_0(\U_G(X,Y))$ of $\U_G(X,Y)$.
\item[{\rm (iii)}] There is a categorical equivalence $\nu_G\colon\Add(\V_G,\Sett)\overset{\simeq}{\longrightarrow}\TamG$.
\end{enumerate}
\end{enumerate}
Here, $\Add(\U_G,\Sett)$ denotes the category of functors from $\U_G$ to $\Sett$, preserving finite products. Similarly for $\Add(\V_G,\Sett)$.
\end{fact}

\begin{cor}\label{CorTX}
By Theorem \ref{MainThm}, there exists an isomorphism of Tambara functors
\[ \nu_G(\V_G(G/G,-))\cong\Omega{[}\mathbf{X}{]} \]
by Yoneda's Lemma.

More generally, if we are given a Tambara functor $T\in\Ob(\TamG)$, then for any $A\in\Ob(\Gs)$, Tambara functor $T_{[A]}=T\underset{\Omega}{\otimes}\nu_G(\V_G(A,-))$ satisfies
\[ \TTamG(T_{[A]},S)\cong S(A) \]
for each $T$-Tambara functor $S$.
Thus $T_{[A]}$ can be regarded as a generalization of $T[\mathbf{X}]$, which represents the evaluation at $A$.
\end{cor}

\begin{rem}
A direct construction of an isomorphism $\Omega[\mathbf{X}]\cong \V_G(G/G,-)$ is also possible. In fact, in the notation of \cite{N_TamMack}, we can construct an isomorphism of semi-Tambara functors $\mathcal{S}(\mathbf{X})\cong\mu_G(\U_G(G/G,-))$, which yields an isomorphism of Tambara functors
\[ \Omega[\mathbf{X}]=\gamma\circ\mathcal{S}(\mathbf{X})\cong\gamma(\mu_G(\U_G(G/G,-)))\cong\nu_G(\U_G(G/G,-)). \]
(We omit the detail. Here, $\gamma\colon\STamG\rightarrow\TamG$ is the left adjoint of the inclusion functor $\TamG\hookrightarrow\STamG$, and $\mathcal{S}\colon\SMackG\rightarrow\STamG$ is the left adjoint of $(\ )^{\mu}\colon\STamG\rightarrow\SMackG$ \ (\cite{Tam},\cite{N_TamMack}).)
\end{rem}

\smallskip

\subsection{A simultaneous generalization of $T[\mathbf{X}]$ and $T[\mathbf{x}]$}

$\ \ $

\smallskip

In the rest, we show the following.
\begin{thm}\label{ThmSimulPoly}
Let $G$ be a finite group, and $T$ be any Tambara functor on $G$. For any $H\le G$, there exists a semi-Mackey functor $\mathfrak{X}_H\in\Ob(\SMackG)$ which satisfies the following.
\begin{itemize}
\item[$(\ast)$] For any $T$-Tambara functor $S$, there exists a natural bijection
\[ \TTamG(T{[}\mathfrak{X}_H{]},S)\cong S(G/H)^{N_G(H)/H}. \]
\end{itemize}
Here, $N_G(H)\le G$ denotes the normalizer of $H$ in $G$. Also, as in Proposition \ref{Prop2}, we have a functor
\[ \mathit{po\ell}_{\mathfrak{X}_H}=\mathcal{F}(-,\mathfrak{X}_H)\colon\TamG\rightarrow\TamG \]
which gives $\mathit{po\ell}_{\mathfrak{X}_H}(T)=T[\mathfrak{X}_H]$ for any $T\in\Ob(\TamG)$.
\end{thm}

Comparing this with Proposition \ref{Prop2}, we have 
\begin{eqnarray*}
T{[}\mathfrak{X}_G{]}&\cong&T{[}\mathbf{X}{]},\\
T{[}\mathfrak{X}_e{]}&\cong&T{[}\mathbf{x}{]},
\end{eqnarray*}
which means Theorem \ref{ThmSimulPoly} generalizes both two situations in Proposition \ref{Prop2}. (More directly, by construction, we have isomorphisms of semi-Mackey functors $\mathfrak{X}_G\cong\mathfrak{A}^{\alpha}$ and $\mathfrak{X}_e\cong\mathcal{L}^G_{\mathbb{N}}$.)

\begin{proof}
Let $H\le G$ be any subgroup, and put $K=N_G(H)/H$. There is a $G$-$K$-biset $U=\mathrm{Indinf}^G_K$, defined to be a set $U=G/H$ equipped with the natural left $G$-action and right $K$-action (\cite{Bouc_Biset}).
Composition by $U$ induces a functor
\begin{eqnarray*}
-\!\circ\! U\colon\SMackG&\rightarrow&\SMackK\\
M&\mapsto&M(U\!\underset{K}{\circ}\! -),
\end{eqnarray*}
which admits a left adjoint $L_U\colon\SMackK\rightarrow\SMackG$ (cf. \cite{N_Biset}, or \cite{Bouc} for the case of Mackey functors).

Remark that we have an isomorphism of $G$-sets $U\!\underset{K}{\circ}\! (K/e)\cong G/H$.

Composing with the functor $\mathcal{L}^K\colon\Mon\rightarrow\SMackK$ in Remark \ref{RemMonAdj}, we obtain a functor
\[ L_U\circ\mathcal{L}^K\colon\Mon\rightarrow\SMackG, \]
which is left adjoint to
\[ \SMackG\overset{-\circ U}{\longrightarrow}\SMackK\overset{\mathit{ev}^K}{\longrightarrow}\Mon. \]
If we put $\mathfrak{X}_H=L_U(\mathcal{L}^K_{\mathbb{N}})$, then we obtain a sequence of natural bijections
\begin{eqnarray*}
\TTamG(T{[}\mathfrak{X}_H{]},S)&\cong&\TamG(\Omega{[}\mathfrak{X}_H{]},S)\\
&\cong& \SMackG(\mathfrak{X}_H,S^{\mu})\\
&=& \SMackG(L_U(\mathcal{L}^K_{\mathbb{N}}),S^{\mu})\\
&\cong& \SMackK(\mathcal{L}^K_{\mathbb{N}},S^{\mu}(U\!\circ\! -))\\
&\cong& \Mon(\mathbb{N},S^{\mu}(U\!\underset{K}{\circ}\! (K/e))^K)\\
&\cong& S^{\mu}(G/H)^{N_G(H)/H}
\end{eqnarray*}
for any $T$-Tambara functor $S$.
\end{proof}

\end{document}